\documentclass[11pt,a4paper]{article}
\usepackage[english]{babel}       \usepackage[latin1]{inputenc}
\usepackage[T1]{fontenc}          \usepackage{a4}
\usepackage[cyr]{aeguill}         \usepackage{epsfig}
\usepackage{amsmath,amsthm,amsfonts,amssymb,mathrsfs}
\usepackage{dsfont}                \usepackage{float}
\usepackage{fancyhdr}              \usepackage{varioref}
\usepackage{color}
\usepackage{url}                   \usepackage{calc}
\usepackage{wrapfig}               \usepackage{graphicx}
\usepackage{enumitem,hyperref}
\usepackage{pifont}       \usepackage{tikz}
\newtheorem{thm}{Theorem} 
\newtheorem{dfn}{Definition}[section]

\newtheorem{prop}{Proposition}[section]

\newtheorem*{conjecture}{Conjecture}
\numberwithin{equation}{section}

\usepackage{authblk}

\makeatletter
\def\namedlabel#1#2{\begingroup
    #2%
    \def\@currentlabel{#2}%
    \phantomsection\label{#1}\endgroup
}
\makeatother

\newcommand{\ensemblenombre}[1]{\mathbb{#1}}

\newcommand{\R}{\ensemblenombre{R}}
\newcommand{\Leb}{\mathcal{L}^d}
\newcommand{\GibbsArea}{\mathcal{G}^{area}}
\newcommand{\Specification}{\mathcal{P}}
\newcommand{\PartitionFunction}{Z^{area}}

\begin{document}

\title{Numerical study for the phase transition of the area-interaction model}

\author{Pierre Houdebert}
\affil{Universität Potsdam, Institut für Mathematik}
\affil{
 \texttt{pierre.houdebert@gmail.com}}
 \date{}
\maketitle

\begin{abstract}
{ 
In this paper we present numerical analysis of the phase transition of the area-interaction model, which is a standard model of Statistical Mechanics.
The theoretical results are based on a recent paper  by Dereudre \& Houdebert 
\cite{dereudre_houdebert_2018_SharpTransitionWR}
which provides a complete phase diagram except on a bounded (implicit) domain.
With our numerical analysis we give an approximative explicit description of this domain.
Furthermore our numerical results confirm the still unproven conjecture stating that non-uniqueness holds if and only if 
$z= \beta$ is large enough, with a value of the threshold obtained from the simulation of $\beta_c \simeq 1.726$.
} 
  \bigskip

\noindent {\it Key words: Gibbs point process, area-interaction, percolation, phase transition, disagreement percolation, random cluster model, Fortuin-Kasteleyn representation, Monte-Carlo algorithm.} 
\end{abstract}

\section{Introduction} \label{section_introduction}
The finite volume area-interaction measure (also called Widom-Rowlinson measure) on a bounded window $\Lambda\subset \R^d$ is defined as a modification of the stationary Poisson Point process of intensity $z$.
Its unnormalized density is given by
$ \exp(-\beta \Leb ( B_1(\omega) ))$
 where $B_1(\omega)$ is the union of the unit balls centred at each point of the configuration $\omega$ and $\Leb$ is the Lebesgue measure on $\R^d$.
 The parameter $\beta\geq 0$ is called \emph{inverse temperature} and the interaction is getting more and more attractive as $\beta$ is large.
The parameter $z$, called \emph{activity}, is related to the intensity of the model.

In  the  infinite  volume  regime  a  global  density  is  senseless and area-interaction measures are defined through equations specifying their conditional laws.
But a solution of these equations can heuristically be seen as the limit of a finite volume area-interaction measure on increasing windows, with an additional boundary condition.
In this paper we are interested in the uniqueness/non-uniqueness (called \emph{phase transition}) of area-interaction measures for given parameters $z$ and $\beta$.
Earlier work of Ruelle \cite{ruelle_1971} proved that in the symmetric case 
$z = \beta$, phase transition occurs when $z = \beta$ is large enough.
A modern proof of this result, based on percolation tools, was done in \cite{chayes_kotecky}.
Until recently almost nothing was proven for the case $z \not = \beta$, and it is conjectured that phase transition (i.e. non-uniqueness) occurs if and \emph{only if} $z=\beta$ large enough.
This conjecture is based on a similar result for the Ising model, see \cite{Velenik}. 
But to the best of our knowledge, no numerical study has been done in order to observe if this conjecture is true.

Recently it was proven in 
\cite{dereudre_houdebert_2018_SharpTransitionWR} that uniqueness of the area-interaction measure is valid for 
$z < \widetilde{z}_c^a(\beta, 1)$, where 
$\beta \mapsto \widetilde{z}_c^a(\beta, 1)$ is a non-decreasing function which is the percolation threshold corresponding to the area-interaction model.
Furthermore 
this function satisfies 
$\widetilde{z}_c^a(\beta, 1) \equiv  \beta$ for $\beta$ large enough.
With some duality property of the model, it provides an almost complete picture of the phase diagram of the model.
But on a bounded region of the parameters ($z, \beta$) it is still not proved whether there is phase transition or not.
We refer to this region as the \emph{unknown region}.
This result is rigorously  stated in Theorem \ref{theo_phase_transition_and_perco} and a sketch of the proof is provided.
Finally the theoretical phase transition diagram is provided in Figure \ref{figure_ph_uniqueness}. 

In the present paper we provide a numerical study of the area-interaction model in order to 
 experimentally plot the curve 
$\beta \mapsto \widetilde{z}_c^a (\beta,1)$ to observe the unknown region;
we validate experimentally the conjecture and find an approximative value of the threshold.
Our simulations are done in dimension $d=2$ using a standard birth and death MCMC algorithm, as presented in \cite{Moller_Waagepetersen_book_2003}.
To implement it we used the so-called \emph{Fortuin-Kasteleyn} representation of the area-interaction model, which provides a construction of the model using the \emph{generalized Continuum Random Cluster Model} (gRCM), which was recently introduced in \cite{houdebert_2019_Potts}.
The definition of the gRCM and the Fortuin-Kasteleyn representation is done at the beginning of Section \ref{section_numerical_study}.
From our numerical study we compute the value of the percolation threshold 
$\widetilde{z}_c^a (\beta,1)$ for several values $\beta$ and provide a plot of the function
$\beta \mapsto \widetilde{z}_c^a (\beta,1)$, see Figure \ref{figure_ph_perco}
and Figure \ref{figure_ph_perco_threshold}.
From this we observe that the unknown region is really small.
Furthermore from our simulation we can validate the conjecture, 
see Figure \ref{figure_ph_intensity_exp},
with an approximative threshold $\beta_c \simeq 1.726$.
This value matches  with numerical studies that have been done in the symmetric case, see \cite{haggstrom_lieshout_moller,
johnson_gould_machta_chayes__MonteCarloStudyWidomRowlinsonFluidUsingClusterMethods}.
\section{Preliminaries} \label{section_preliminaries}
Let us consider the state space $\R ^d $.
Let $\Omega$ be the set of locally finite configurations $\omega$ on $\R^d$.
This means that $\# (\omega \cap \Lambda)< \infty$ 
for every bounded Borel set $\Lambda$ of $\R^d$, with $\#\omega$ being the cardinality of the configuration $\omega$.
We write $\omega_{\Lambda}$ as a shorthand for $\omega \cap \Lambda$.
To a configuration $\omega \in \Omega$ we associate the germ-grain structure
$
B_r(\omega) :=
\underset{x \in \omega}{\bigcup} B_r(x),
$
where $B_r(x)$ is the closed ball centred at $x$ with radius $r>0$.
Let $\pi^{z}$ be the distribution on $\Omega$ of the standard homogeneous Poisson point process with intensity $z >0$. 
For $\Lambda \subset \R^d$ bounded, we denote by $\pi^{z}_\Lambda$  the restriction of  
$\pi^{z}$  on $\Lambda$.
\subsection{Area-interaction measures}
The area-interaction measures - also called Widom-Rowlinson measures - are defined through the standard Gibbs DLR formalism prescribing the conditional probabilities.
For a bounded $\Lambda \subset \R^d$, we define the \emph{$\Lambda$-Hamiltonian}
$$
H_{\Lambda} (\omega)
:=
\Leb ( \ B_1(\omega_{\Lambda}) \setminus B_1(\omega_{\Lambda^c}) \ )
$$
with $\Leb$ the standard $d$-dimensional Lebesgue measure.
The \emph{area specification} on a bounded $\Lambda \subseteq \R^d$ with boundary condition $\omega_{\Lambda^c}$ is defined by
\begin{align*}
\Specification_{\Lambda, \omega_{\Lambda^c}}^{z,\beta}(d \omega'_\Lambda)
: =
\frac{
e^{-\beta H_{\Lambda}
(\omega'_{\Lambda}  \omega_{\Lambda^c}) }}
{\PartitionFunction_{z,\beta ,\Lambda,\omega_{\Lambda^c}} } \pi^{z}_{\Lambda}(d \omega'_\Lambda)
\end{align*}
with the standard \emph{partition function}
$$
\PartitionFunction_{z,\beta ,\Lambda,\omega_{\Lambda^c}}
: =
\int_{\Omega}
\ e^{-\beta H_{\Lambda}
(\omega'_{\Lambda}  \omega_{\Lambda^c}) }
\pi^{z}_{\Lambda}(d\omega'_{\Lambda})
$$
which is always non-degenerate 
(i.e. $0<\PartitionFunction_{z,\beta ,\Lambda,\omega_{\Lambda^c}}<+\infty$).
Let us point out that for $\beta=0$, we have
$\Specification_{\Lambda, \omega_{\Lambda^c}}^{z,\beta}
\equiv 
\pi^{z}_{\Lambda}$.
\begin{dfn}
\label{def_area_model}
A probability measure $P$ on $\Omega$ is an \emph{area-interaction measure} of activity $z$ and inverse temperature $\beta$,  written 
$P \in \GibbsArea_{z, \beta}$, if for every bounded Borel set $\Lambda \subset \R^d$ and  every bounded measurable function $f$,
\begin{align}
\label{eq_dlr_area_model}
\int_{\Omega} f \ d P 
=
\int_{\Omega} \int_{\Omega}
f(\omega'_{\Lambda}  \omega_{\Lambda^c} )
\Specification^{z, \beta}_{\Lambda, \omega_{\Lambda^c}}(d \omega'_\Lambda)
P (d \omega ).
\end{align}
The equations \eqref{eq_dlr_area_model}, for all bounded $\Lambda$, are called \emph{DLR equations}, after Dobrushin, Lanford and Ruelle.
Those equations prescribe the conditional probabilities of a Gibbs measure.
\end{dfn}
Heuristically a solution of the DLR equations can be seen as the limit of
$\Specification_{\Lambda_n, \omega_{\Lambda_n^c}}^{z,\beta}$ for an increasing sequence $\Lambda_n$ and some boundary condition $\omega$.
There is \emph{phase transition}, i.e. non-uniqueness, if the limit depends on the boundary condition $\omega$.
\subsection{Percolation}
The theory of percolation studies the connectivity in random structures and is a crucial tool to prove phase transition of the area-interaction measure.
\begin{dfn}
Let $r>0$.
A configuration $\omega$ is said to $r$-percolate if the germ-grain structure 
$B_r(\omega)$ has at least one unbounded connected component.
Furthermore
a probability measure $P$ on $\Omega$ is said to \emph{$r$-percolate} (resp. \emph{do not percolate}) if 
$P(\{ \omega  \text{ $r$-percolates} \})=1$ 
(resp. $P(\{ \omega \text{ $r$-percolates} \})=0$).
\end{dfn}
Thanks to standard monotonicity arguments applied to the Gibbs specification,
we have:
\begin{prop}
\label{propo_threshold_area}
For all $\beta>0$ and $r>0$, there exists 
$0<\widetilde{z}_c^a(\beta,r) <\infty$ such that
\begin{itemize}
\item
for all $z < \widetilde{z}_c^a(\beta,r)$, 
every area-interaction measure $P \in \GibbsArea_{z, \beta}$
almost never $r$-percolates, i.e
$
P(\{ \omega \ r\text{-percolates}\})
= 0;
$
\item
for all $z > \widetilde{z}_c^a(\beta,r)$, every area-interaction measures $P \in \GibbsArea_{z, \beta}$ almost surely $r$-percolates, i.e
$
P(\{ \omega \ r\text{-percolates}\})
= 1.
$
\end{itemize}
\end{prop}
A proof of this result is provided in 
\cite[Prop. 2.7]{dereudre_houdebert_2018_SharpTransitionWR}.
In the general case the only information known about 
$\widetilde{z}_c^a(\beta,r)$ is the following bound coming from stochastic domination:
\begin{align*}
\widetilde{z}_c^a(0,r)
\ \leq \
\widetilde{z}_c^a(\beta,r)
\ \leq \
\widetilde{z}_c^a(0,r)\exp (\beta v_d),
\end{align*}
where $\widetilde{z}_c^a(0,r)$ is the percolation threshold of the Poisson Boolean model of constant radii $r$.
Experimental studies showed that in dimension $d=2$ we have
$\widetilde{z}_c^a(0,r) \simeq 0.359072 \cdot r^{-2}$, see
\cite{Quintanilla__Torquato_measurement_percolation_threshold_2000}.
To the best of our knowledge there exists no approximation for cases 
$\beta \not =0$.
\section{Theoretical results} \label{section_results}
The first fundamental question in Gibbs point process theory is the existence of at least one probability measure satisfying the DLR equations \eqref{eq_dlr_area_model}.
This is an interesting and non trivial question treated for several kinds of interactions, see e.g.
\cite{dereudre_drouilhet_georgii,
dereudre_houdebert,
dereudre_vasseur_2019_existence,
roelly_zass_2019_existence}.
Since the area-interaction process has a finite range interaction, existence is long proved, see \cite{ruelle_livre_1969}.
\begin{prop}
For all $z, \beta \geq 0$, the set $\GibbsArea_{z, \beta}$ of area-interaction measures is non-empty.
\end{prop}
The second question concerning the area-interaction process is its uniqueness/ non-uniqueness, known as \emph{phase transition}.
It is conjectured that non-uniqueness happens if and only if $z=\beta$ is large enough.
\begin{conjecture}
There exists $0<\beta_c<\infty$ such that phase transition occurs for the area-interaction model if and only if $z = \beta > \beta_c$.
\end{conjecture}
This conjecture is motivated by a similar result proved for the lattice Ising model, see for instance \cite[Th. 3.28 \& Th. 3.46]{Velenik}.
Although this conjecture is still open, the following theorem based on our recent work \cite{dereudre_houdebert_2018_SharpTransitionWR} provides an almost complete picture of the phase diagram,
see Figure \ref{figure_ph_uniqueness}.
\begin{thm}
\label{theo_phase_transition_and_perco}
${}$
\begin{enumerate}
\item
There exists $\beta_1<\infty$ such that for $z=\beta>\beta_1$, there is non-uniqueness of the area-interaction measure;
\item  for $z, \beta$ such that 
$z < \widetilde{z}_c^a(\beta,1)$, there is uniqueness of the area-interaction measure.
Using a duality property, we obtain the same for 
$\beta<\widetilde{z}_c^a(z, 1)$;
\item
there exists $\beta_2 \in [\beta_1,\infty[$ such that 
$\widetilde{z}_c^a(\beta,1) = \beta$ for $\beta>\beta_2$.
This partially proves the conjecture.
\end{enumerate}
\end{thm}
\begin{figure}[t]
\centering
\begin{tikzpicture}[scale=0.5]
\draw [>=latex,->] (0,0) -- (10,0) node[right] {$z$};
\draw [>=latex,->] (0,0) -- (0,10) node[above] {$\beta$};
\draw [dashed] (0,0) -- (10,10) node[above right] {$z = \beta$};
\draw (5,5) node {$\bullet$};
\draw [dashed] (5,5) -- (5,0);
\draw (5,0) node {} node[below] {$\beta_1$};
\draw [red] (5,5) -- (10,10);
\draw [>=latex,<-] (7.5,7.5) -- (11,9.5); 
\node[draw,align=center] at (14,8.5) {Non uniqueness \\ symmetric case };
\draw [line width=1pt] (2,0) arc (180:120.05:8)--(7.5,7.725);
\draw [line width=1pt] (7.5,7.725) arc (300:315:5)--(10,10);

\draw (8,8) node {$\bullet$};
\draw [dashed] (8,8) -- (8,0);
\draw (8,0) node {} node[below] {$\beta_2$};

\node[align=center] at (-3,4) 
{$\beta \mapsto \widetilde{z}_c^a(\beta,1)$};
\draw [>=latex,<-] (2.7,3.5) -- (-0.5,4);

\fill [blue, opacity=0.60] (0,0)--(2,0)-- (2,0) arc (180:120.05:8)--
(7.5,7.725)--(7.5,7.725) arc (300:315:5)-- (10,10)--(0,10);

\fill [blue, opacity=0.60](0,0)-- (0,2)-- (0,2) arc (270:329.5:8)--
(7.725,7.5)--(7.725,7.5) arc (150:135:5)-- (10,10)--(10,0);

\node[draw,align=center] at (5,11.5) {Uniqueness};
\draw [>=latex,<-] (3,7.5) -- (5,10.5);
\draw [>=latex,<-] (7,3) -- (5,10.5);
\draw (4,3) node[above] {\bf \Huge ?} ;
\end{tikzpicture}
\caption{Theoretical uniqueness/non-uniqueness regimes for the area-interaction measures with parameters $z,\beta$.}
\label{figure_ph_uniqueness}
\end{figure}
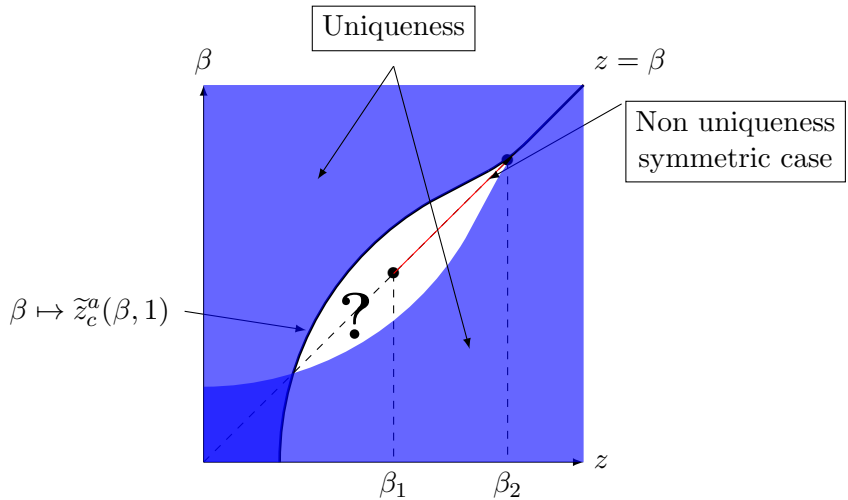
\begin{proof}[Sketch of the proof]
The first point was initially proved in 1971 by Ruelle  (see \cite{ruelle_1971}) using Peierls argument.
Then a modern proof based on a Fortuin-Kasteleyn representation satisfied by the area-interaction model was provided in \cite{chayes_kotecky}.
As a bi-product of their construction, the authors proved that 
$\widetilde{z}_c^a(\beta,1/2) =\beta$ for $\beta>\beta_1$.
The second point is a generalisation of the disagreement percolation construction introduced in \cite{Hofer-temmel_Houdebert_2018}.
The idea of disagreement percolation is to compare the Gibbs specification 
$\Specification^{z, \beta}_{\Lambda, \omega_{\Lambda^c}}$ with the same parameters 
$z, \beta$ but with two different boundary conditions 
$\omega^{1}_{\Lambda^c},\omega^{2}_{\Lambda^c}$, using percolation ideas.
The duality property is a consequence of the representation using the Widom-Rowlinson model, and is stated in 
\cite[Proposition 2.5]{dereudre_houdebert_2018_SharpTransitionWR}.
The third point is proved in
\cite[Section 4.4.1]{dereudre_houdebert_2018_SharpTransitionWR} using the Fortuin-Kasteleyn representation introduced in \cite{chayes_kotecky} and an elegant stochastic domination argument which, using the fact that as soon as 
$\widetilde{z}_c^a(\beta,1/2) =\beta$ for $\beta$ large, one gets that 
$\widetilde{z}_c^a(\beta,1) =\beta$ for $\beta$ even larger.
\end{proof}
\section{Numerical study of the phase diagram}
\label{section_numerical_study}
In this section we will use numerical approximation, in dimension $d=2$, in order to
\begin{enumerate}
\item experimentally plot the curve 
$\beta \mapsto \widetilde{z}_c^a (\beta,1)$ to see the region which is not covered by the Theorem \ref{theo_phase_transition_and_perco}.
\item validate experimentally the conjecture and find an approximative value of the threshold $\beta_1$;
\end{enumerate}

To do our numerical study we will use a birth and death MCMC algorithm to sample the area-interaction process.
The general algorithm  we used for birth and death MCMC can be found in 
\cite{Moller_Waagepetersen_book_2003}.
One could have considered using an {\it exact simulation} technique, as implemented in \cite{haggstrom_lieshout_moller}, but the computation of the percolation threshold requires to sample the model in a large window, which would be extremely time consuming using exact simulation techniques.
We chose to sample the area-interaction model from the Fortuin-Kasteleyn representation using the gRCM, introduced in
\cite{houdebert_2019_Potts}.
\begin{dfn}
On a bounded window $\Lambda$, the generalized Continuum Random Cluster Modem with activity parameter $\rho \geq 0$ and $\alpha_1,\alpha_2 \geq 0$ such that $\alpha_1 + \alpha_2=1$ is defined as
$$
P_{\Lambda}^{gRCM}(d \omega)
\sim
\prod_{C}
\left( 
\alpha_1^{\# C} + \alpha_2^{\# C}
\right)
\pi_{\Lambda}^{\rho}(d\omega),
$$
where the product is over the clusters of $B_{1/2}(\omega)$.
\end{dfn}
\begin{prop}
Considering a configuration $\omega \sim P_{\Lambda}^{gRCM}$ and removing each cluster $C$ of $B_{1/2}(\omega)$ with probability
$$
{\alpha_2^{\# C}}/({\alpha_1^{\# C} + \alpha_2^{\# C}}),
$$
one obtains a configuration sampled from the area-interaction measure 
$$
P_{\Lambda}^{area}(d \omega)
\sim
\exp (- \beta |B_1(\omega) \cap \Lambda|)
\pi_{\Lambda}^{z}(d\omega)
$$
 with parameters 
$z= \alpha_1 \rho  $ and $\beta = \alpha_2 \rho$.
\end{prop}
This representation gives a good feeling for the third point of Theorem \ref{theo_phase_transition_and_perco}.
Indeed a large cluster will be removed (resp. kept)  with very high probability when 
$\alpha_1<\alpha_2$ (resp. $\alpha_1>\alpha_2$).
But percolation is highly dependent of the status of the large clusters.

In our numerical study we sample  area-interaction measures in a window $\Lambda=[0,100]^2$, using a MCMC algorithm sampling the gRCM and then thinning the configuration according to the previous proposition. We observe the experimental intensity and whereas the center of the window is connected to its boundary in $B_1(\omega)$.
For each pair of parameters $(z, \beta)$, we sampled the model 
{\bf 1000} times in order to obtain an experimental intensity and an experimental probability of percolation.
The C++ code used is accessible on GitHub\footnote{\url{https://github.com/PierreHoudebert/area_perco_multithread}}.
\begin{figure}[h!]
\begin{center}
\includegraphics[width=6cm]{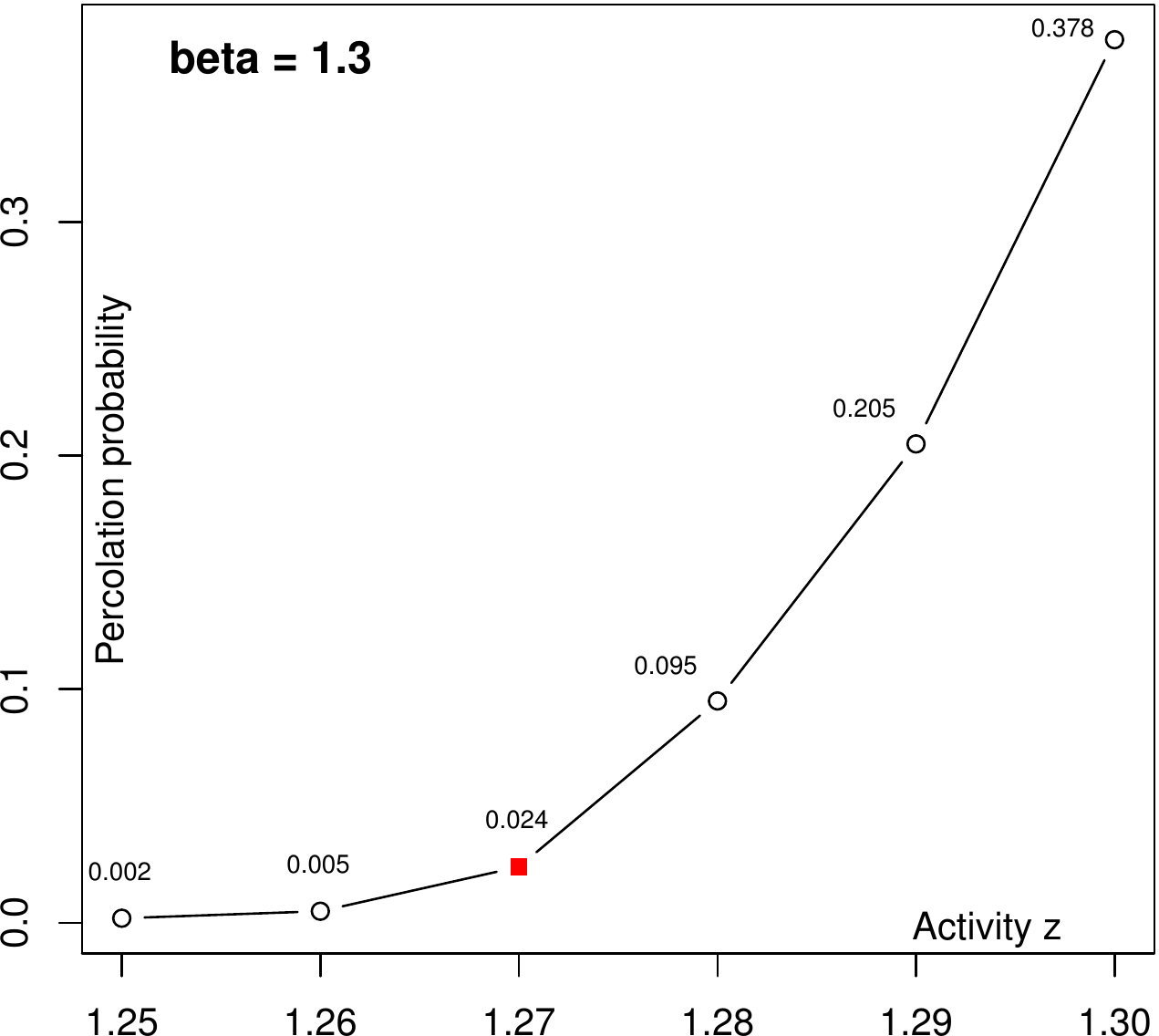}
\includegraphics[width=6cm]{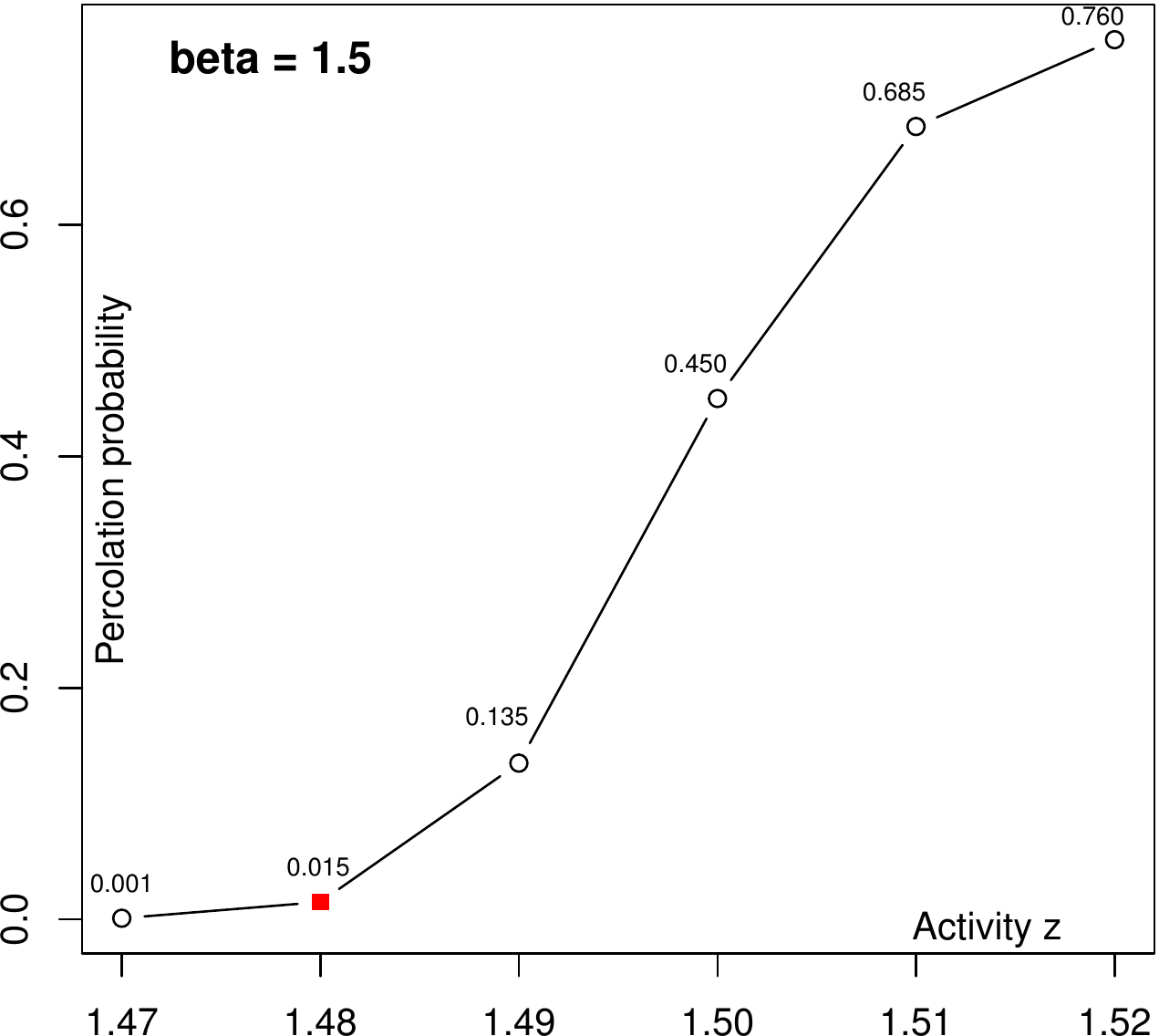}

\vspace{0.2cm}

\includegraphics[width=6cm]{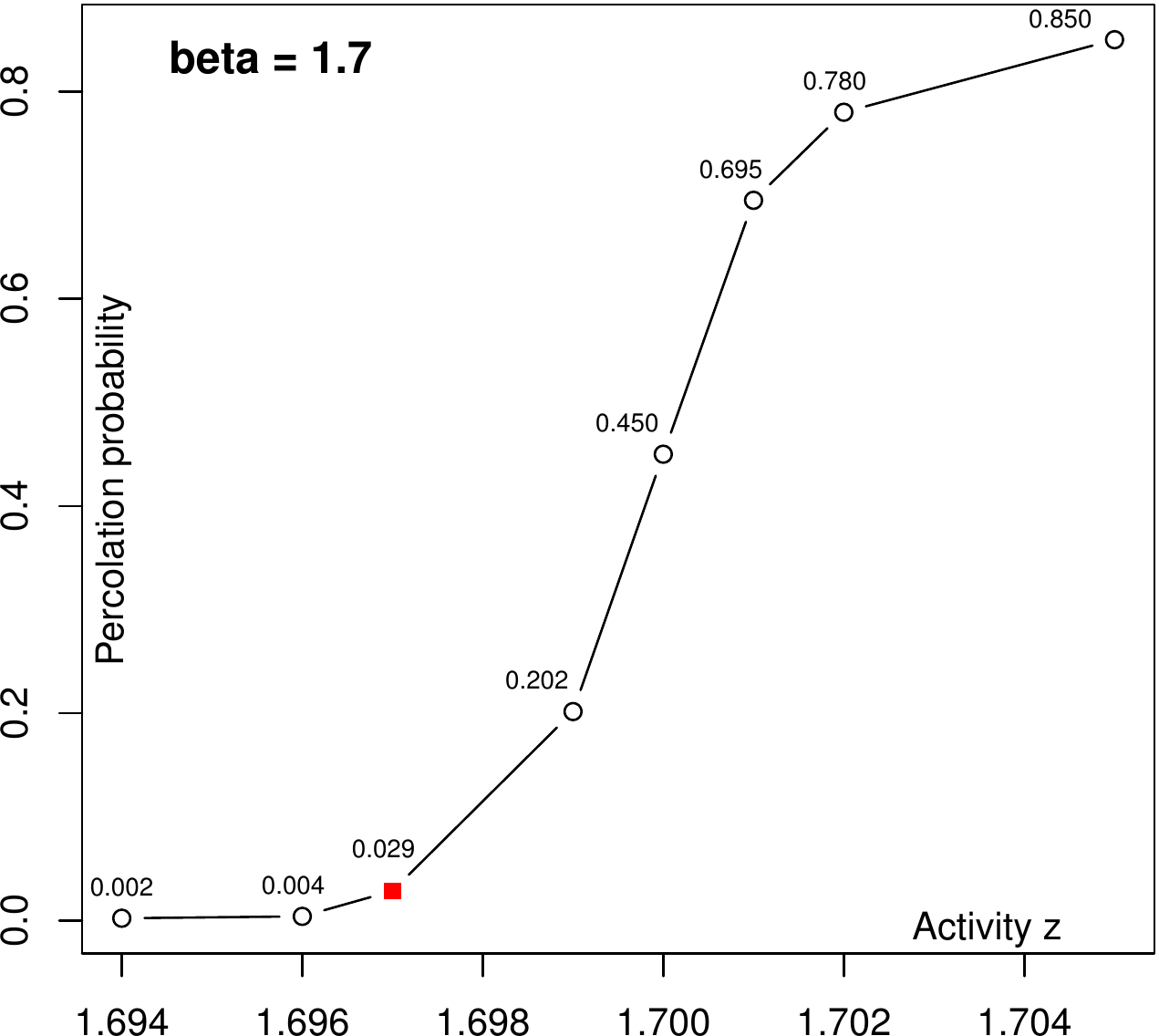}
\includegraphics[width=6cm]{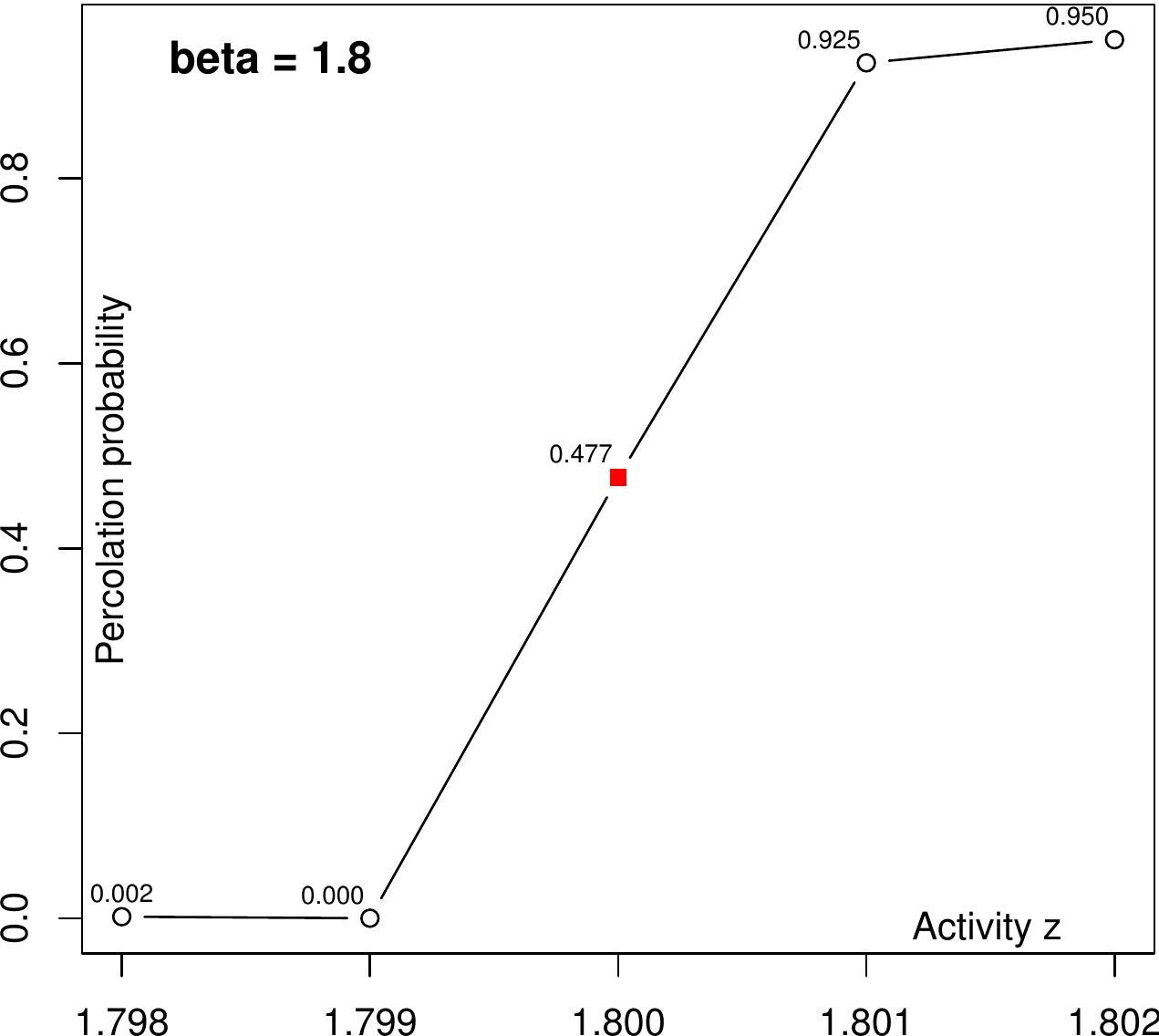}
\end{center}
\caption{Percolation probability as a function of the activity $z$. The value corresponding to the percolation threshold $\widetilde{z}_c^a (\beta,1)$ 
is the plain square.}
\label{figure_ph_perco}
\end{figure}
\begin{figure}[h!]
\begin{center}
\begin{minipage}{8cm}
\includegraphics[width=5cm]{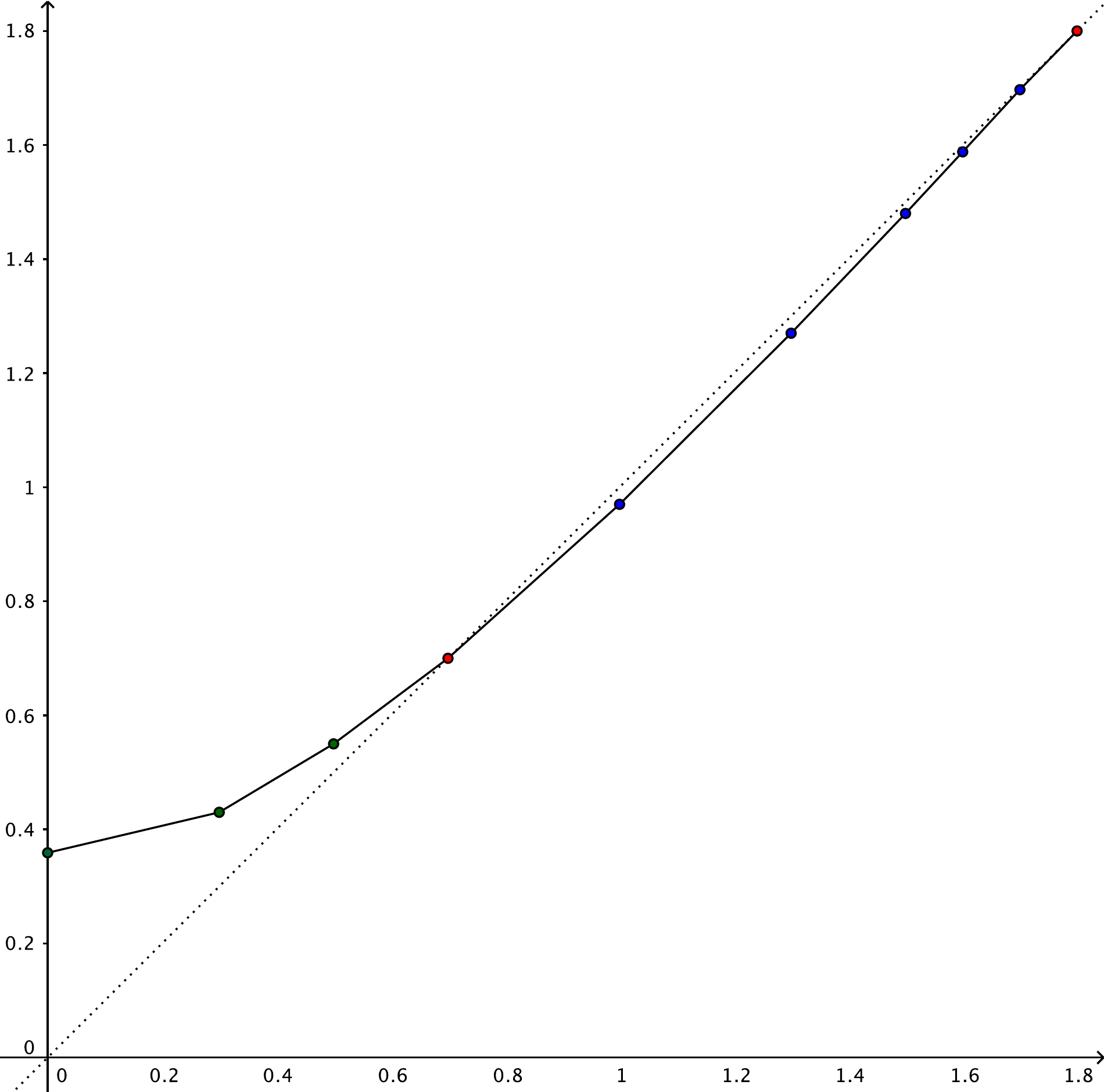}
\end{minipage}
\begin{tabular}{|c|c|}
\hline
$\beta$ & $\widetilde{z}_c^a (\beta,1)$
\\ \hline
$0.3$ &  $0.43$
\\
$0.5$ &  $0.55$
\\
$0.7$ &  $0.7$
\\
$1$    &  $0.97$
\\
$1.3$ &  $1.27$
\\
$1.4$ &  $1.37$
\\
$1.5$ &  $1.48$
\\
$1.6$ &  $1.588$
\\
$1.7$ &  $1.697$
\\
$1.72$ &  $1.718$
\\
$1.725$ &  $1.723$
\\
$1.8$ &  $1.8$
\\
$2$    & $2$
\\ \hline
\end{tabular}
\caption{Obtained values of the percolation threshold, and corresponding approximation curve}
\label{figure_ph_perco_threshold}
\end{center}
\end{figure}
In order to determine $\widetilde{z}_c^a (\beta,1)$, the percolation threshold 
of the model for a given $\beta$, we arbitrarily decide that the percolation threshold is the first observed value of $z$ such that the probability of the center of the box $[0,100]^2$ to be connected to the boundary is larger than $0.01$.
Indeed, theoretically, for the infinite volume model the threshold is the first value of $z$ such that the probability of the origin (or any given point) belongs to the infinite connected component is positive.
But for the finite volume model this probability is always positive.

For some values of $\beta$ the percolation probability, as a function of $z$, is displayed in Figure \ref{figure_ph_perco}.
The obtained values of the percolation threshold 
$\widetilde{z}_c^a (\beta,1)$
and the corresponding graph is displayed in Figure \ref{figure_ph_perco_threshold}.

Considering the conjecture.
It is known that phase transition occurs at a given pair $(z_0,\beta_0)$ if the intensity, as a function of $z$ for $\beta_0$ fixed,  is discontinuous in 
$z_0$.
This is proven for the Ising model, see for instance \cite{Velenik}, and could be proven similarly for the area-interaction model.

\begin{figure}[h]
\centering
\includegraphics[width=13cm]{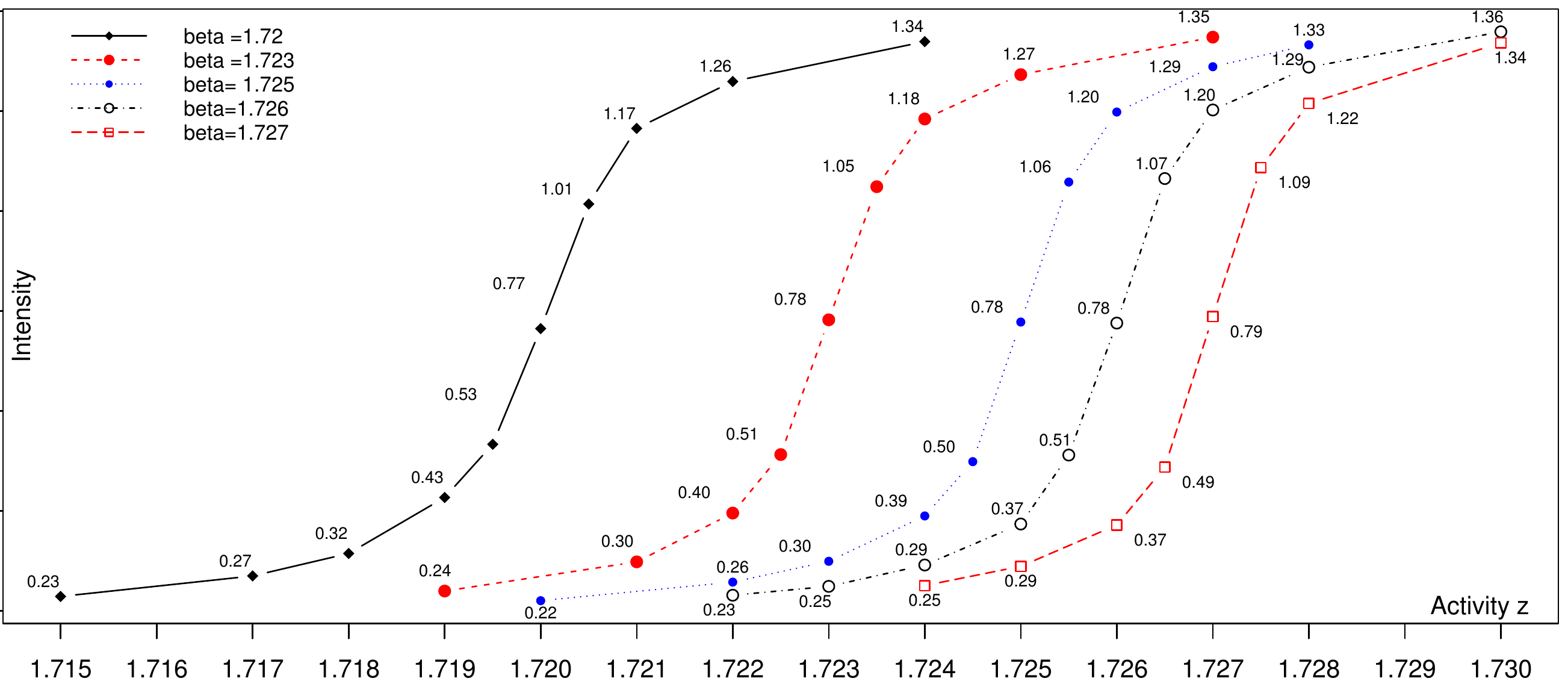}
\caption{Experimental intensity as a function of $z$, for several values of $\beta$. }
\label{figure_ph_intensity_exp}
\end{figure}

As before for a given pair of parameters $(z, \beta)$, we sample the model 1000 times on the bounded window $\Lambda=[0,100]^2$, and observed the experimental intensity obtained.
For given values of $\beta$, we provide the graph of the intensity as a function of $z$ in Figure \ref{figure_ph_intensity_exp}.
We observe from Figure \ref{figure_ph_intensity_exp}
that the experimental intensity is
indeed discontinuous only for $z =\beta$ larger than some value $\beta_c$.
This amounts to saying that phase transition occurs only for $z = \beta$, as conjectured.
However from the figure it is not clear what is the exact value of the threshold $\beta_c$, but it seems to be approximately $\beta_c \simeq 1.726$, which is coherent with the value obtained in \cite{johnson_gould_machta_chayes__MonteCarloStudyWidomRowlinsonFluidUsingClusterMethods} where the authors  considered only the symmetric case.

\vspace{1cm}

{\it Acknowledgement:} 
This  work was  supported  in part by the ANR  project PPPP (ANR-16-CE40-0016) and by Deutsche Forschungsgemeinschaft (DFG) - SFB1294/1 - 318763901.
\bibliographystyle{plain}
\bibliography{biblio}
\end{document}